\documentclass[12pt]{article}
\def\const{\mathrm{const}}

\def\Cn{{\overline{C}_n}}
\def\C{{\mathbf{C}}}
\author{Alexandre Eremenko\thanks{Supported by NSF grant DMS-1361836.}}
\title{Zeros and coefficients}
\begin{document}
\maketitle
\begin{abstract} Two theorems  on the asymptotic
distribution of zeros of sequences of
analytic functions are proved. First one relates the asymptotic behavior
of zeros to the asymptotic behavior of coefficients. Second theorem establishes
a relation between the asymptotic behaviors of zeros of a function
and zeros of derivative.
\end{abstract}

Let $f_n$ be a sequence of analytic functions in a region
$D$ in the complex plane. The case of entire functions, $D=\C$
will be the most important.
Suppose that there exist positive numbers $V_n\to\infty$,
such that the subharmonic functions
\begin{equation}\label{1}
\frac{1}{V_n}\log|f_n|\quad \mbox{are uniformly bounded from above}
\end{equation}
on every compact subset of $D$, and for some $z_0\in D$ the sequence
\begin{equation}\label{2}
\frac{1}{V_n}\log|f_n(z_0)|\quad\mbox{is bounded from below.}
\end{equation}
Under these conditions, one can select a subsequence such that
for this subsequence
\begin{equation}\label{3}
\frac{1}{V_n}\log|f_n|\to u,
\end{equation}
where $u$ is a subharmonic function in $D$, $u\neq-\infty$.
This convergence can be understood in various senses, for example,
quasi-everywhere, which means at every point except a set of points of
zero logarithmic capacity, and also in $D'$ (Schwartz distributions),
or in $L^p_{\mathrm{loc}}$, $1\leq p<\infty$.

The corresponding Riesz measures converge weakly. The limit measure,
which is $(2\pi)^{-1}\Delta u$, describes the limit distribution
of zeros of $f_n$.

For all these facts we refer to \cite{H1} or \cite{H2}.

This setting is useful in many cases when one is interested in asymptotic
distribution of zeros of sequences of analytic functions.

Let us mention several situations to which our results apply.
\vspace{.1in}

1. Let $f_n$ be monic polynomials of degree $n$, and $V_n=n$.
We denote by $\mu_n$ the so-called empirical measures, which are the
Riesz
measures of $(1/n)\log|f_n|$. The measure $\mu_n$ is discrete, with finitely
many atoms, with an atom of mass $m$ at every zero of $f_n$ of multiplicity $m$.
So the $\mu_n$ are probability measures in the plane. Therefore,
one can always choose a subsequence for which $\mu_n\to\mu$, where $\mu$
is some measure in the plane. The sequence of potentials will converge:
$$\frac{1}{n}\log|f_n|\to u=\int_{|\zeta|\leq 1}\log|z-\zeta| d\mu(\zeta)+
\int_{|\zeta|>1}\log\left|1-z/\zeta\right| d\mu(\zeta)+C.$$
\vspace{.1in}

2. Let $f$ be an entire function, $f(0)\neq 0$, and 
\begin{equation}\label{M}
M(r,f)=\max\{|f(z)|:|z|=r\}
\end{equation}
its maximum modulus. Suppose that for some sequence $r_n\to\infty$ we
have
$$\log M(2r_n,f)=O(\log M(r_n,f)),\quad n\to\infty.$$
Then the functions $f_n(z)=f(r_nz)$ satisfy (\ref{1}) and (\ref{2})
with $V_n=\log M(r_n)$ and $z_0=0$.
\vspace{.1in}

3. The order of an entire function is defined by
$$\rho=\limsup_{r\to\infty}\frac{\log\log M(r,f)}{\log r}.$$
An entire function of finite order is said to be of normal type if
$$\limsup_{r\to\infty}\frac{\log M(r,f)}{r^\rho}<\infty.$$
If $f$ is of finite order and normal type, $f(0)\neq 0$,
and $r_n$ is any sequence tending
to infinity, then one can take $f_n(z)=f(r_nz)$, and $V_n=r_n^{\rho}$,
and the properties (\ref{1}) and (\ref{2}) (with $z_0=0$) will hold.

Such sequences of subharmonic functions were applied to the asymptotic study
of zeros in many papers. Of the early ones, we mention \cite{Delange}.
See also
\cite{Azarin}, \cite{BKW}, \cite{Brolin}, \cite{Eremenko}, \cite{KZ},
\cite{Ronkin}, \cite{Sodin}
and literature in these books and papers.
  
The purpose of this note is to state and prove two general facts about this
convergence, and give several applications.
There is no claim of novelty of the main results or their proofs, our only goal
is to bring these general facts to the attention of a wider audience,
as on our opinion they deserve.

The result of section 3 was found in the 1990s jointly with M. Sodin.
This section is independent on section 1.

I thank M. Sodin, A. Rashkovskii and D. Novikov for useful discussions.

\section{Maximum modulus and coefficients}

In this section the functions $f_n$ are defined and analytic in
a disc $D=\{ z:|z|<R\}$ where $0<R\leq\infty$, and we assume (\ref{1}),
(\ref{2}) and (\ref{3}).

Let $M(r,f_n)$ be the maximum modulus, and $B(r)=B(r,u)=\max_{|z|=r}u(z)$,
then we have $\log M(r,f_n)/V_n\to B(r),\; n\to\infty$.
This convergence is uniform on every interval $[0,r_0]$ with $r_0<R$.
Let 
$$\Phi(t)=B(e^t), \quad -\infty<t<\log R.$$
Suppose that
\begin{equation}\label{fk}
f_n(z)=\sum_{k=0}^\infty a_{n,k}z^k,
\end{equation}
and let 
$$\psi_n(x)=\sup_{\mbox{convex}\;\psi}
\left\{ \psi:\psi\left(\frac{k}{V_n}\right)\leq
-\frac{\log|a_{n,k}|}{V_n}\right\}.$$
Here the sup is taken over all convex functions $\psi$ on the real line.
To check that the class is not empty, we use the assumption that our functions
$f_n$ are analytic in some disc centered at $0$ which implies that
$\log|a_{n,k}|\leq C_nk$, so $\psi(x)=-C_nx$ is in the class.

It is convenient to consider all convex functions as defined
on the whole real line and taking values in $(-\infty,+\infty]$.
A convex function with finite values on a half-line or
on an interval (closed, open, semi-open)
can be always extended
to the whole real line by setting it equal to $+\infty$ at every point
outside this interval.

The set where a convex function is finite is always convex (but can be
neither open nor closed).
Then pointwise $\sup$ of any family of convex functions is convex,
and the limit of decreasing sequence of convex functions is convex
if finite at least at one point. Thus our functions $\psi_n$ are convex.

The set where a convex function is discontinuous can consist
of at most two points, the endpoints of the interval where it is finite.
Lower semi-continuous regularization is still convex
and differs from the original function at at most $2$ points.
We always assume our convex functions to be lower semi-continuous, making
this regularization when necessary.
\vspace{.1in}

\noindent
{\bf Theorem 1.} {\em Suppose that (\ref{3}) holds.
Then $\psi_n\to\phi>-\infty$, and $\phi$ is
the Legendre transform of $\Phi$,
$$\phi(x)=\sup_t(tx-\Phi(t)).$$
Conversely, if $\liminf_{n\to\infty}\psi_n=\phi>-\infty,$
then (\ref{1}) holds
and one can choose a subsequence for which 
(\ref{3}) holds, and $\psi_n\to\phi$, and $\Phi$ and $\phi$
are related by the Legendre transform.}
\vspace{.1in}

{\em Proof.} First we recall that Legendre's transform is defined
for convex functions of the real
line. 
The Legendre transform is always a convex, lower semi-continuous function.
The transform is order-reversing: $f\geq g$ implies $L(f)\leq L(g)$.
Finally, the Legendre transform is an involution, $L\circ L={\mathrm{id}}$.

We have by Cauchy inequality
$$|a_{n,k}|r^k\leq M(r,f_n),$$
Applying log, dividing by $V_n$, and writing $\log r=t$, we rewrite this as
$$\frac{-\log|a_{n,k}|}{V_n}\geq \sup_t\left(\frac{k}{V_n}t-\frac{1}{V_n}
\log M(e^t,f_n)\right).$$
We set set $x=k/V_n$, and as the RHS is convex in $x$,
we obtain that
$$\psi_n(x)\geq \sup_t\left(xt-\log M(e^t,f_n)\right)/V_n.$$
We pass to the limit when $n\to\infty$ while $x$
is fixed,
and obtain
\begin{equation}\label{one}
\phi(x)\geq\sup_t(xt-\Phi(t)).
\end{equation}

Now we need an estimate in the opposite direction.
Let us first give the estimate for the important special case when
$f_n$ are polynomials of degree $n$, and $V_n=n$.
In this case our argument is very simple.
We start with the trivial inequality
\begin{equation}\label{ca}
n\max_k|a_{n,k}|r^k\geq M(r,f_n).
\end{equation}
Taking logs, and substituting $\log r=t,$ we obtain
$$\max_k\left(kt+\log|a_{n,k}|\right)\geq \log M(e^t,f_n)-\log n.$$
Dividing by $V_n=n$ and passing to the limit as $n\to\infty$
with our usual arrangement
that $k/n=x$, we obtain
\begin{equation}\label{6}
\max_x(xt-\phi(x))\geq \Phi(t).
\end{equation}
Now, Legendre's transform is an involution, 
and 
order
reversing. So our two inequalities imply the equality: $\phi$ and $\Phi$
are Legendre transforms of each other.

Returning to the general case, we need a substitute for (\ref{ca}).
This is provided by rudimentary Wiman--Valiron theory, which we recall.
Let 
$$f(z)=\sum_{k=0}^\infty a_kz^k$$
be a power series with radius of convergence $R\leq+\infty$. For every
$r\in[0,R)$, we introduce the maximal term
$$m(r,f)=\max_{k}|a_k|r^k,$$
and the central index $\nu(r)$ as the largest integer for which
$$m(r,f)=|a_{\nu(r)}|r^{\nu(r)}.$$
For $0<r<r_1<R$, we have
\begin{equation}
\label{A}
\nu(r)\log\frac{r_1}{r}\leq\log m(r_1)-\log m(r).
\end{equation}

For our functions $f_n$ as in (\ref{fk}), notice the relation
\begin{equation}\label{mr}
\frac{1}{V_n}\log m(e^t,f_n)=\max_x(xt-\psi_n(x))=L(\psi_n)(x).
\end{equation}

\noindent
{\bf Lemma 1.} (Valiron). {\em For every $r$ and $r_1$ such that
$$0\leq r<r_1<R,$$
we have
\begin{equation}\label{B}
M(r,f)\leq m(r_1,f)\left(\nu(r)+\frac{r_1}{r_1-r}\right).
\end{equation}}
\vspace{.1in}

{\em Proof.}
\begin{eqnarray*}
M(r,f)&\leq&\sum_{k=0}^\infty|a_k|r^k\\
&\leq& m(r,f)\nu(r)+\sum_{k=\nu(r)+1}^\infty
|a_k|r_1^k\left(\frac{r}{r_1}\right)^k\\
&\leq&
m(r_1,f)\nu(r)+m(r_1,f)\frac{r_1}{r_1-r}.
\end{eqnarray*}

This is a sufficient substitute for (\ref{ca}) for our purpose.

Combining (\ref{A}) and (\ref{B}), we obtain
$$M(r,f_n)
\leq m(r_1,f_n)\left(\frac{\log m(r_1,f_n)-\log m(r,f_n)}{\log r_1-\log r}
+\log r_1-\log(r_1-r_1)\right).$$
Taking logarithms,
dividing by $V_n$ and passing to the limit as $n\to\infty$,
we take into account that
$$\lim_{n\to\infty}(\log m(r,f_n))/V_n=L(\phi)(\log r),$$
and obtain that
$$\Phi(t)\leq L(\phi)(t'),$$
which holds for every $t<t'$. As both sides are continuous, we conclude
that
$$\Phi\leq L(\phi)$$ which is (\ref{6}).
This completes the proof.
\vspace{.1in}

Theorem 1 gives an asymptotic connection between the maximum modulus and
coefficients.

There is one interesting case when the function $B(r)$ completely determines
the limit zero distribution $\mu$.
\vspace{.1in}

\noindent
{\bf Proposition.} {\em
Suppose that the function $z\mapsto B(|z|,u)$ is piecewise-harmonic
in $\{ z:|z|<R\}$ in the following sense:
there exists a closed nowhere dense set $E\subset(0,r)$ such that for every
component $I$ of $(0,R)\backslash E$, we have
$$B(r,u)=a_I\log r+b_I,\quad r\in I,$$
and $B(r,u)=\const$ for sufficiently small $r$.
Then $u(z)=B(|z|,u),\; |z|<R$.}
\vspace{.1in}

For example, if
$B(r,u)=\log^+r$ then $u(z)=\log^+|z|$.
We conclude that $\mu$ is the uniform distribution
on the unit circle if and only if $B(r,u)=\log^+|z|$.
The corresponding functions are
$\Phi(t)=t^+,\; -\infty<t<+\infty$,
and $\phi(x)=0,\; 0\leq x\leq 1$. 
\vspace{.1in}

We obtain  necessary and sufficient conditions, in terms of the coefficients,
for the zeros to be uniformly distributed on the unit circle.
For simplicity, we state it only for sequences of polynomials.
\vspace{.1in}

\noindent
{\bf Corollary}. {\em Let 
$$f_n=\sum_{k=0}^na_{n,k}z^k$$
be a sequence of polynomials, and $\mu_n$
the empirical measures of $f_n$.
Then $\mu_n$ converge weakly to the uniform
probability measure on the unit circle, as $n\to\infty$, if
and only if the following condition holds.
For every $\epsilon>0$ there exists $n_\epsilon$, such that for every
$n>n_\epsilon$
we have
\begin{equation}\label{cond}
\max_{k\in[0,n]}\log|a_{n,k}|-\max_{k\in[0,\epsilon n]\cup[(1-\epsilon)n,n]}
\log|a_{n,k}|
\leq\epsilon n.
\end{equation}
}
\vspace{.1in}

An equivalent result was recently obtained by Fernandez \cite{F}.

For example, condition (\ref{cond}) holds if all coefficients are
uniformly bounded from above, while $|a_{n,0}|$ and $|a_{n,n}|$
are bounded from below.

To prove the Proposition we notice that on any complementary
interval $I$, $U(z):=a_I\log|z|+b_I$ is a harmonic majorant
of $u$, and $u(z)=U(z)$ for some point inside the ring
$\{ z:|z|\in I\}$.
So my the maximum principle $u=U$.

\section{Applications of the previous results}

\noindent
{\bf Example 1.} (Jentzsch)
Let $f$ be a power series with radius of convergence $1$
and $f_{n_k}$ a sequence of partial sums such that $|a_{n_k}|^{1/n_k}\to 1$.
Then $f_{n_k}$ satisfy the assumptions of the Proposition, and
we obtain a theorem of Jentzsch that zeros of $f_{n_k}$ are
asymptotically uniformly distributed on the unit circle.
\vspace{.1in}

\noindent
{\bf Example 2.} Asymptotic distribution of zeros of a Ruelle zeta-function
\cite{ELS}. 
This zeta function is related to the dynamical system $z\mapsto p_c(z)=z^2+c$,
where $c$ is a real parameter.
The zeta function is
$$F_c(z)=1+\sum_{k=1}^\infty\frac{z^k}{p_c(0)\ldots p_c^{*n}(0)},$$
where $p_c^{*n}$ is the $n$-th iterate of $p_c$. For $c<-2$,
we have a real one-parametric family of entire functions. Zeros of these
functions are eigenvalues of a Ruelle operator. Let $V(c)$ be the smallest
positive integer $n$, for which $p_c^{*(n+1)}/p_c^{*n}\geq 36$.
By analyzing Taylor's coefficients of $F_c$, it was shown in \cite{ELS}
that
$$\frac{1}{V(c)}\log M(F_c,r)\to \log^+(r/2),\quad c\to -2,\quad r>0.$$
Therefore, the Proposition with $R=\infty$ implies that the limit distribution of zeros
is the uniform measure on the circle $|z|=2.$ 

Similar arguments apply to the entire function
$$H_a(z)=\frac{1}{a}\sum_{n=0}^\infty a^{2^n}z^n,\quad 0<a<1,$$
studied by Hardy \cite{H}. As $a\to 1-,$ the limit distribution of zeros
is the uniform measure on the unit circle.
\vspace{.1in}

\noindent
{\bf Example 3.} 
The results about this example were first obtained by Alan Sokal and
Alexander Scott (2005, unpublished). 

Consider the polynomials
$$\Cn(y)=\sum_{k=n-1}^{n(n-1)/2}c_{n,k}(y-1)^k,$$
where $c_{n,k}$ is the number of simple connected graphs
with $k$ edges on $n$ labeled vertices.
Tutte \cite{Tutte} discovered the following identity for
the the exponential generating function of $\Cn$:
\begin{equation}
\label{11}
\sum_{n=1}^\infty\Cn(y)x^n/n!=\log F(x,y),
\end{equation}
where
$$F(x,y)=\sum_{n=0}^\infty y^{n(n-1)/2}x^n/n!.$$
This series $F$ represents an entire function of $x$ for $|y|\leq 1$.
It is of order zero when $|y|<1$, and thus it has roots.
If $x_0(y)$ is the root of the smallest modulus, then
the series in the left hand side of (\ref{11}) 
is convergent for $|x|<|x_0(y)|$.

{}From the formula for the radius of convergence, we obtain
$$\limsup_{n\to\infty}n^{-1}\log|\Cn(y)/n!|=-\log|x_0(y)|,\quad |y|<1.$$
In it not difficult to prove that actually the limit
exists
\begin{equation}\label{rco2}
\lim_{n\to\infty}n^{-1}\log|\Cn(y)/n!|=-\log|x_0(y)|\quad\mbox{for}\quad
|y|<1.
\end{equation}

This relation justifies the need to study the function $y\mapsto x_0(y)$.
For a survey of known results and many conjectures about this function
we refer to \cite{Sok}.

{}From the definition of $\Cn$ we infer that
these are monic polynomials of degree $d_n=n(n-1)/2$. So
\begin{eqnarray*}
u_n(y)&:=&d_n^{-1}\log|\Cn(y)/n!|
=
d_n^{-1}\sum_{k=0}^{n(n-1)/2}\log|y-y_k|+c_n\\
&=&
\int\log|y-w|d\mu_n(w)+c_n,
\end{eqnarray*}
where $\mu_n$ is the empirical measures of
$\Cn$, and $c_n\sim (\log n!)/d_n$
are constants tending to zero. From the combinatorial definition
one can compute the values
$\Cn(0)=(-1)^{n-1}(n-1)!$. Thus 
\begin{equation}\label{geq}
\lim_{n\to\infty}u_n(0)=0.
\end{equation} 
Equation (\ref{rco2}) implies convergence $u_n\to u$,
where $u$ is the potential of the limit measure $\mu$ of
empirical measures.

As the plane is not compact, the limit measure $\mu$ can {\'a} priori
have smaller total mass than $1$,
but we will see that this does not happen in our situation:
the zeros of $\Cn$ do not escape
to $\infty$.

We have $u_n\to u$, where
\begin{equation}
\label{pot}
u(y)=\int\log|y-w|d\mu(w),\quad \mu(\C)\leq 1.
\end{equation}
Equation (\ref{rco2}) implies that 
\begin{equation}
\label{leq0}
u(y)\leq 0,\quad \mbox{for}\quad |y|<1.
\end{equation}
On the other hand, (\ref{geq}) implies
that $u(0)\geq 0$,
so, by the Maximum Principle,
\begin{equation}
\label{zero0}
u(y)=0\quad\mbox{for}\quad |y|<1.
\end{equation}
By a well-known property of subharmonic functions,
this implies
\begin{equation}\label{zero}
u(y)=0\quad\mbox{for}\quad |y|\leq 1.
\end{equation}
On the other hand, from (\ref{pot}) follows that $u(y)\leq \log|y|+o(1)$
as $y\to\infty$, and together with
(\ref{zero}) this implies that for every $\epsilon>0$, the function 
$(1+\epsilon)\log|y|$ is a harmonic majorant of $u$ in the annulus
$1<|y|<\infty$. So we have
\begin{equation}
\label{Leq}
u(y)\leq \log^+|y|,\quad\mbox{for all}\quad y\in \C.
\end{equation}
We already know that equality holds in the unit disc. To show that
it holds everywhere, we produce a lower estimate as follows.
{}From Cauchy's inequality
$$2^{d_n}\leq M(2,\Cn):=\max_{|y|=2}|\Cn(y)|,$$
and it follows that $M(2,u)\geq \log 2.$
Thus we have equality in (\ref{Leq}) at some point $y_1$, with $|y_1|=2$.
By the Maximum Principle, applied to the ring $1<|y|<\infty$
we have equality in (\ref{Leq}) everywhere, that is
\begin{equation}\label{logplus}
u(y)=\log^+|y|.
\end{equation}
This implies that $\mu=(2\pi)^{-1}\Delta u$
is the uniform distribution on the unit circle, in particular $\mu$
is a probability measure in $\C$ as advertised above. Moreover,
$\mu_n$ tends to the uniform distribution on the unit circle.

\section{Asymptotic behavior of derivatives}

In this section, $D$ is an arbitrary region in the complex plane,
and $f_n$ analytic functions in $D$.
Suppose that the sequence of subharmonic functions $(\log|f_n|)/n$
converges:
$$\frac{1}{n}\log|f_n|\to u.$$
Then one can choose a subsequence such that
$$\frac{1}{n}\log|f^\prime_n|\to v.$$
We fix a sequence for which both limits hold.
\vspace{.1in}

\noindent
{\bf Theorem 2.} {\em Under the stated conditions, $v\leq u$, and
if at some point $v(z_0)<u(z_0)$, then there is a neighborhood
$V$ of $z_0$ such that $u$ is constant on $V$.}
\vspace{.1in}

{\em Proof.} For every subharmonic function $w$ and $\epsilon>0$ we define
$$w_\epsilon (z)=\max_{\zeta:|\zeta-z|\leq\epsilon}w(\zeta).$$
Then $w_\epsilon$ is a decreasing sequence of subharmonic functions,
and for every $z$, $w(z)\leq w_\epsilon(z) \to w(z)$ as $\epsilon\to 0$.

Choose an arbitrary $z_0\in D$. Cauchy's estimate gives
$$|f_n^\prime(z)|\leq\epsilon^{-1}
|f_n|_{2\epsilon}(z_0),\quad|z-z_0|<\epsilon.$$
Taking logs, dividing by $n$ and passing to the limit, we conclude that
$$v_\epsilon(z_0)\leq u_{2\epsilon}(z_0).$$
This holds for every sufficiently small $\epsilon>0$, so passing to the limit
when $\epsilon\to 0$, we obtain $v(z_0)\leq u(z_0)$.

To prove the second statement of the Theorem, we fix $z_0\in D$, and
suppose that 
\begin{equation}\label{0}
u(z_0)-v(z_0)=5\delta>0.
\end{equation}

By upper semi-continuity of $v_0$, there exists a disc $V_{2r}$ of radius
$2r$ centered 
at $z_0$, such that
$$v(z)<v(z_0)+\delta,\quad z\in V_{2r}.$$
this implies that there exists $n_0$ such that for
$n>n_0$ we have
\begin{equation}\label{111}
\log|f_n^\prime(z)|\leq n(v(z_0)+2\delta),\quad z\in V_r.
\end{equation}

Consider the set 
$$E=\{ z\in V_r:u(z)\geq u(z_0)-\delta\}.$$
This set is non-empty because it contains $z_0$, and it has positive area.

According to a theorem of Azarin \cite{Azarin}, the set
$$\{ z\in V_r:|u(z)-(\log|f_n|)/n|>\delta\}$$
can be covered by a countable set of discs with the sum of
 radii tending to $0$
as $n\to\infty$. This means that for $n$ sufficiently large, there exist
points $z_n\in E$ such that
\begin{equation}\label{22}
|f_n(z_n)|\geq n(u(z_0)-2\delta).
\end{equation}
Now for $z\in V_{r}$ we have in view of (\ref{111}):
$$\left||f_n(z)|-|f_n(z_n)|\right|\leq|f_n(z)-f_n(z_n)|\leq
2r e^{n(v(z_0)+2\delta)}.$$
This gives in view of (\ref{22}) and (\ref{0}) that
$$\log|f_n(z)|=\log|f_n(z_n)|+\log\left(1+2re^{-n\delta}\right)=
\log|f_n(z_n)|+o(1),$$
which implies that $u$ is constant in $V_r$.
\vspace{.1in}

{\bf Example 4.} Let $P_n$ be polynomials of degree $n$ and $\mu_n$
their empirical measures. Suppose that $\mu_n\to\mu$ weakly,
where $\mu$ is a Borel measure in the plane.
\vspace{.1in}

{\em If the complement of the closed support of $\mu$
 does not have bounded components, then
the empirical measures $\mu^\prime$ describing the distribution of
critical points of $P_n$ converge to the same measure $\mu$.}

Purdue University

West Lafayette, IN 47907

eremenko@math.purdue.edu
\end{document}